\documentclass[12pt]{article}

\usepackage{amsfonts}
\usepackage{amssymb}
\usepackage{amsmath}
\usepackage{amsthm}
\usepackage{epsfig}

\title{An analog to Deuring's criterion for good reduction of elliptic curves
\thanks{Work supported in part by the European Community's Human 
Potential Programme
under Contract HPRN-CT-2000-00114, GTEM}}

\author{Claus Lehr}

\newcommand{\PK}{{\mathbb P}^{1}_{K}}

\newtheorem{defi}{Definition}[section]
\newtheorem{proposition}[defi]{Proposition}
\newtheorem{theorem}[defi]{Theorem}
\newtheorem{lemma}[defi]{Lemma}

\newtheorem{remark}[defi]{Remark}

\newcommand{\B}{\beta}
\newcommand{\C}{\gamma}
\newcommand{\E}{\epsilon}
\newcommand{\U}{x_2}
\newcommand{\ed}{\hfill$\square$}

\newcommand{\AP}{{\mathbb P}^{1}_{R}}

\newcommand{\llll}{\lambda}

\newcommand{\Pk}{{\mathbb P}^{1}_{k}}

\newcommand{\la}{\longrightarrow}

\begin{document}
\maketitle

\begin{abstract}
In this paper we study the reduction of $p$-cyclic covers of the $p$-adic
line ramified at exactly four points. For $p=2$ these covers are
elliptic curves and Deuring has given a criterion for when such a
curve has good reduction. Here we consider the case of $p>2$ and 
completely determine the stable model of the cover.
In particular we obtain a finite extension $R'$ of $R$ necessary
for the stable reduction to be defined.
No additional conditions are imposed on the geometry of the branch locus
and thus this work can be viewed as a first step towards understanding
the situation where branch points coalesce.
\end{abstract}

\section{Introduction}
Let $R$ be a complete mixed characteristic $(0,p)$ discrete valuation
ring. We denote by $K$ the field of fractions of $R$ and by $k$ the 
residue field of $R$ which we assume to be algebraically closed.
The valuation given on $K$ by the ring $R$ will be called $v$.
For $X \la \PK$ a $p$-cyclic cover we are interested in the stable
model of $X$. In previous work (cf. \cite{L1}) we have considered this 
problem under
certain conditions on the geometry of the branch locus $B$ of the
cover. In particular we required that $B$ has equidistant geometry, which is
to say that $\PK$ has a smooth $R$-model such that $B$ extends to a 
relative \'etale divisor over $R$. This condition was first introduced
in a paper by Raynaud (cf. \cite{R1}) and
it is also a key hypothesis in \cite{M} generalizing the results in 
\cite{L1}. 
Here we allow any $B$ consisting of 
four rational points and determine the stable reduction $X_k$ of the
resulting cover $X$.
The assumption of $B$ to be equidistant implies in particular that the 
stable reduction has no vanishing cycles - a result essentially due to
Raynaud (cf. \cite{R1}). Knowing this the computations are then reduced
to finding the components of positive genus in the stable reduction. 
Under the conditions imposed here 
this will no longer be true, in particular $X$ can be a Mumford curve.
The situation of a general branch locus $B$ has been dealt with in
\cite{L2}. There we gave explicit methods to determine the number of 
cycles. Beyond this for such a general branch locus it is still not
known how to obtain a fine description of the stable reduction of $X$.\\

\noindent{\bf Acknowledgments}  
I would like to thank David Harbater for suggesting the example given
in section \ref{exo}. He also pointed out problems in the notation of
a previous version. \\

\section{The main result} \label{exres}
We keep the notation of the introduction and assume in addition that $R$
contains a primitive $p$-th root of unity $\zeta$ as well as 
$\tau=(-p)^{p/(p-1)}$. We denote by $\pi$ a uniformizer in $R$. 
Let $X \la \PK$ be a $p$-cyclic cover. Suppose the cover is ramified at
exactly four $K$-rational points. For a suitable coordinate $x_0$ on $\PK$
we can assume the branch locus to be $\{0,1,\infty, \llll\}$ with
$\llll \in K-\{0,1\}$. Similarly making an appropriate change of coordinates
we can assume $\llll \in R-\{0,1\}$ such that the reduction 
$\bar{\llll}\not=1$.
Therefore the cover is given birationally by the equation 
$z_0^p=x_0(x_0-1)^{\B}(x_0-\llll)^{\C}$ with positive integers $\B,\C < p$ such that 
$(1+\B+\C,p)=1$. 

The following theorem
completely describes the structure of the stable reduction of $X$. Its 
strength lies in the fact that we do not only determine the combinatorics
of the special fiber (which can be done with less work)
but also give equations for the stable model $X_{R'}$
over an explicitly determined discrete valuation ring $R' \supseteq R$.
In particular this yields equations for the irreducible components of the 
special fiber. We write $$j(\llll)= 
p^{-2p/3(p-1)}\left( \llll^2(\B+1)^2-2\llll(\B+\C+1-\B\C)+(\C+1)^2\right).$$

\begin{theorem} \label{deuring}
With the notation introduced above, and for $p>3$,
there are three different possibilities for the special fiber 
of the stable model $X_k$:
\begin{itemize}
\item[1)] $X_k$ consists of a single irreducible component, i.e.\ 
$X$ has potentially
good reduction iff one of the following two conditions is satisfied:

a) $v(\llll)=0$ and  $v(j(\llll)) \ge 0$.

In this case $X_k$ has $p$-rank $0$.

b) $\C+1=p$ and $v(\llll)=v(\tau^2)$.

In this case $X_k$ has $p$-rank $p-1$.

\item[2)] $X_k$ has exactly two irreducible components, both of genus zero and
intersecting each other in $p$ distinct points iff $\C+1=p$ and 
$v(\llll)>v(\tau^2)$. In this case $X$ is a Mumford curve.

\item[3)] In all remaining cases $X_k$ has two irreducible components, each
of genus $(p-1)/2$, intersecting in precisely one point.
\end{itemize}
\end{theorem}

\begin{figure}[h]
\begin{center}
\input{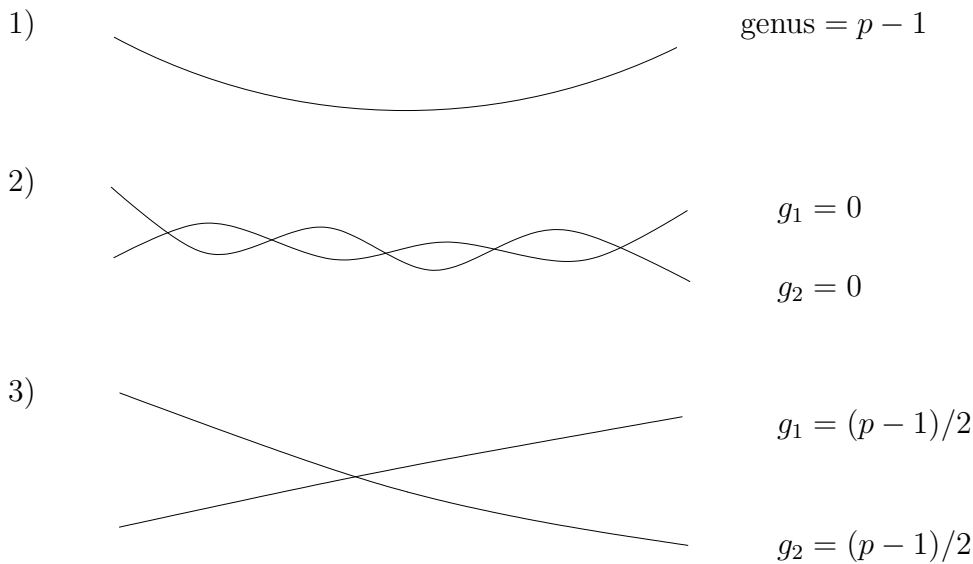}
\caption{The special fiber in the three cases of Theorem \ref{deuring}}
\label{3fibers}
\end{center}
\end{figure}

\begin{remark} \rm In case 1) of Theorem \ref{deuring} $X_k$ is a smooth curve and a $p$-cyclic cover
of $\Pk$. In the sub-case a) the branch locus of this cover has exactly one point while in b) it
consists of two points.
\end{remark}

\begin{defi} \rm For the purposes of this paper we call the different types of degeneration
in Theorem \ref{deuring} type 1a), 1b), 2) and 3).
\end{defi}

\subsection{Semi-Stable Models} \label{mod}
To be precise we recall in this section the definitions of the various
$R$-models considered.

\begin{defi} \rm An $R$-scheme $Y_R$ is called
{\it semi-stable} if it is proper, flat of relative dimension $1$
over $R$ and the following conditions hold:

\noindent i) Its generic fiber $Y_K$ is smooth over $K$.

\noindent ii) Its special fiber $Y_k$ is reduced and
has at most ordinary double points.

\noindent The scheme $Y_R$ is also called semi-stable model for $Y_K$.
\end{defi}

\begin{defi} \rm A semi-stable curve $Y_R$ is called {\it stable}
if every irreducible
component of $Y_k$ that is isomorphic to $\Pk$ meets the other components
in at least $3$ points.
\end{defi}

\begin{remark} \rm For any smooth $K$-curve of genus $\ge 2$ there exists a
stable model, after possibly passing to a finite extension of $K$
(cf.\ \cite{DM}).
\end{remark}

\section{Proof of the main theorem}

In order to show Theorem \ref{deuring} we will proof Theorem 
\ref{deuring}' below from which Theorem \ref{deuring} is easily deduced.
\vspace{.5cm}

{\noindent{\bf Theorem \ref{deuring}'} \it
Keeping the above notation we distinguish two cases:
\begin{itemize}
\item[1)] Assume $v(\llll)=0$. 

a) $v(j(\llll)) \ge 0$ iff $X$ has potentially good reduction of type 1a).

b) $v(j(\llll)) < 0$ iff $X_k$ consists of two irreducible components,
each of genus $(p-1)/2$ intersecting in precisely one point, i.e.\ reduction of type 3).

\item[2)] Assume $v(\llll)>0$.

a) If $\C+1 \not=p$ then $X_k$ has reduction of type 3).


b1) If $\C+1=p$ and $v(\llll)=v(\tau^2)$, then $X$ has potentially good 
reduction of type 1b).

b2) If $\C+1=p$ and $v(\llll)>v(\tau^2)$, then $X_k$ has reduction of type 2).

b3) If $\C+1=p$ and $v(\llll)<v(\tau^2)$, then $X_k$ has reduction of type 3).

\end{itemize}}

The main tool in proving Theorem \ref{deuring} is the following result on
the degeneration of ${\bf \mu}_p$-torsors, which we recall from \cite{L1}.
To state it we introduce the following notation. If $v$ is the valuation 
defined by $R$ on $K$ we extend it to the ring $R[x]$ in the following way:
\[\mbox{For} \quad \sum^m_{i=0} a_ix^i \in R[x] \quad \mbox{define} \quad
v(\sum^m_{i=0} a_ix^i)=\mbox{min}\{v(a_i) | 0 \le i \le m\}. \]

Let $C \la \PK$ be a $p$-cyclic cover given birationally by $y^p=f(x)$ with
$f(x) \in R[x]$ monic. The equation determines an $R$-model $\mathcal C$ for $C$,
an affine patch of which is obtained by normalizing $R[x]$ in the function
field $K(C)$. We further assume that ${\mathcal C}_k={\mathcal C} 
\otimes_R k$ 
is reduced. Under these
assumptions we have the following result characterizing ${\mathcal C}_k$.

\begin{proposition}\label{ap}
Choose $h(x) \in R[x]$ such that $w:=v(h(x)^p-f(x)) \in {\mathbb Z}$ 
is maximal. 

\noindent 1) $w>v(\tau)$:
Then ${\mathcal C}_k$ has $p$ irreducible components each of which is
isomorphic to $\Pk$. The map ${\mathcal C}_k \la \Pk$
is an isomorphism when restricted to each irreducible component.

\noindent 2) $w=v(\tau)$:
The Artin-Schreier equation 
\begin{equation} \label{fi}
T^p-T+\big((h(x)^p-f(x))/\tau\big)^-/\bar{h}(x)^p=0
\end{equation}
is 
irreducible over the field $k(x)$ and ${\mathcal C}_k \la \Pk$ is given
birationally by this equation.

\noindent 3) $w<v(\tau)$: 
Then ${\mathcal C}_k \la \Pk$ is purely inseparable.
\end{proposition} \ed \\

\begin{remark} \label{app} \rm In order to apply Proposition \ref{ap} it is not always necessary to
know the exact value of $w$.
For instance if one can find an $h(x)$ such that $v(h(x)^p-f(x))=v(\tau)$
and equation \eqref{fi} is irreducible it follows that $w=v(\tau)$.

Also if for a certain $h(x)$ with $v(\pi^e)=v(h(x)^p-f(x))<v(\tau)$ the
polynomial $((h(x)^p-f(x))/\pi^e)^- \notin k[x]^p$ then $w<v(\tau)$.
Further in this situation the cover ${\mathcal C}_k \la \Pk$ is given 
birationally by $T^p=((h(x)^p-f(x))/\pi^e)^-$.

Both of these assertions are immediate consequences of the proof of the
proposition.

\end{remark}

We add an easy lemma that will be useful later.

\begin{lemma} \label{hh}
Consider the situation of Proposition \ref{ap} case 2). We write $N=sp$
for the degree of $f(x)$ and assume that
$$f(x)=\sum_{i=1}^N a_Nx^N \quad \mbox{\rm with } v(a_N) \ge v(\tau^{(N-i)/p}).$$
Then $v(x^N-f(x))=v(\tau)$ and hence ${\mathcal C}_k \la \Pk$ is given by
$$T^p-T+\big((x^N-f(x))/\tau\big)^-/x^N=0.$$
\end{lemma}

The following classical result will be used at various points in the paper.
A proof can be found in \cite{St} , III, 7.8.

\begin{proposition} \label{stich} Let $k$ be an algebraically closed field of
characteristic $p>0$ and $f(x) \in k[x]$ a polynomial of
degree $m$ such that the Artin-Schreier equation 
$y^p-y=f(x)$ defines a field extension $L/k(x)$, and $(m,p)=1$. 
Then the genus of $L$ is $(m-1)(p-1)/2$.
\end{proposition} \ed \\

\noindent{\it Proof of Theorem \ref{deuring}'.}
Some constructions in this proof will require to pass from $R$ to a finite
extension $R'$. We will not always point this out but note here that in
the results used each of these extensions is given explicitly.  

1)a): We apply \cite{L1} Theorem 1 with $f(x_0)=x_0(x_0-1)^{\B}(x_0-\llll)^{\C}$,
$m=3$, $n=1+\B+\C$, $b=\tau^{1/3}$ 
and proceed as in the proof of loc.\ cit.\ Corollary\ 2.
We compute $$f'(x_0)=(x_0-1)^{\B-1}(x_0-\llll)^{\C-1}\left( (\B+\C+1)x_0^2
-x_0(\B\llll+\llll+\C+1) +\llll\right) .$$ 
So $g(x_0)=x_0^2-x_0(\llll(\B+1)+\C+1)/(\B+\C+1)+\llll/(\B+\C+1)$ and
$g'(x_0)=2x_0-(\llll(\B+1)+\C+1)/(\B+\C+1)$.
Assuming potentially good reduction 
condition a) in loc.\ cit. Theorem 1 reads:
\begin{equation}\label{pp}
g(x_0) \equiv (x_0-d)^2=x_0^2-2dx_0+d^2 \quad (\mbox{\rm mod } b)
\end{equation} 
On the other hand we have the above expression for $g(x_0)$ and hence
we set $d=(\llll(\B+1)+\C+1)/2(\B+\C+1)$.
Notice that condition \eqref{pp} determines $d$ modulo $b$ and hence, if
there exists a choice of $d$ inducing good reduction, so will our choice of
$d$ above.
Now we evaluate the conditions in part
b) of loc.\ cit. Theorem 1.
With respect to the coordinate $x_1=b/(x_0-d)$ we get the following equation
for the generic fiber:
\begin{equation}\label{blow}
z_1^p=x_1^N+b\frac{f'(d)}{f(d)}x_1^{N-1}+b^2\frac{f^{(2)}(d)}{2f(d)}x_1^{N-2}+
b^3 \frac{f^{(3)}(d)}{3!f(d)}x_1^{N-3}+ \dots
\end{equation}
where $N \in \{p,2p\}$ depends on $n$.
Notice that in the reduction $\bar{f}'(x_0)$ of the polynomial $f'(x_0)$ to $k[x_0]$
the zero $\bar{d}$ has multiplicity $2$ and therefore $f^{(3)}(d)$ is a unit
in $R$. Further $f(d)$ is a unit in $R$.
Now we apply Proposition \ref{ap}. Still assuming we have potentially good reduction
the cover given by \eqref{blow} has to fall into case 2) of the proposition and the genus of the special
fiber is $p-1$. Let $s$ be such that $N=sp$. Then $s$ is the degree of 
$h(x_1)$ and since $p \ge 5$ 
we may assume $h(x_1)=x_1^s$ using Lemma \ref{hh}.
Here $s \in \{1,2\}$ depending on the degree
of $f(x_0)$.
We conclude that potentially good reduction implies
$$v(f'(d))\ge v(b^2) \quad \mbox{\rm and} \quad v(f^{(2)}(d))\ge v(b).$$
Now $v(f'(d))=v(g(d))$ so $v(g(d)) \ge v(b^2)$. Finally
$$g(d)=\frac{-1}{4(\B+\C+1)}(\llll^2(\B+1)^2-2\llll(\B+\C+1-\B\C)+(\C+1)^2)$$
and $4(\B+\C+1)$ is a unit. This shows that the condition on $j(\llll)$ is
necessary. 

To show sufficiency define $d,b$ as before and  
observe that $j(\llll) \ge 0$ implies
$\bar{g}(\bar{d})=0$ so $f(d)$ is a unit. Also by construction $g'(d)=0$.
We conclude that $\bar{d}$ is a zero of multiplicity $2$ in $\bar{f}'(x_0)$
hence $f^{(3)}(d)$ is a unit in $R$.
Define $F(x_0)$ by $f'(x_0)=F(x_0)g(x_0)$. Then, using $v(j(\llll))\ge 0$,
we get $v(f'(d)) \ge v(b^2)$ and the 
identity $f^{(2)}(x_0)=F'(x_0)g(x_0)+F(x_0)g'(x_0)$ shows $v(f^{(2)}(d)) \ge v(b^2)$. 
Now the above equation \eqref{blow} yields a smooth model using 
Proposition \ref{ap} as indicated in Remark \ref{app}.

1)b): It follows from $v(\llll)=0$ that the stable reduction of $X$ has
a tree-like special fiber. Further the sum of the geometric genera over all 
the components in this tree is equal to the genus $g$ of the generic fiber $X$
(cf.\ \cite{L1} Theorem 2). By the Riemann-Hurwitz formula we get
$g=p-1$. As case 1)a) is an 
'if and only if' statement in the present case we will have at least two 
components of strictly
positive genus. On the other hand by Artin-Schreier Theory
the genus of such a component is a 
multiple of $(p-1)/2$ (cf.\ Proposition \ref{stich}).
We conclude that the special fiber of the stable model has two irreducible
components of genus $(p-1)/2$ intersecting in exactly one point. 
The equations of those components have been given in the proof of loc.\ cit.\
Theorem 2.

\begin{figure}[h]
\begin{center}
\input{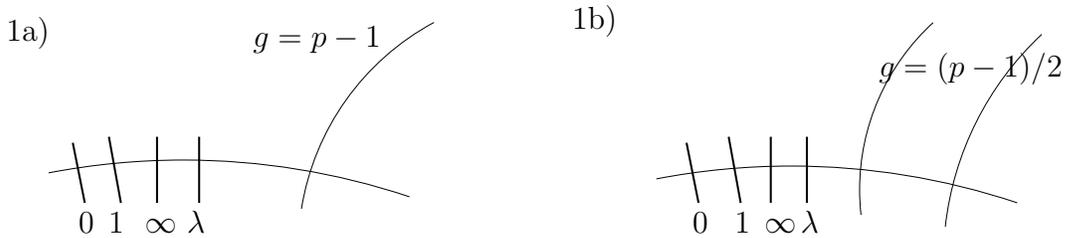}
\caption{Typical picture of the fibers in Theorem \ref{deuring}' 1a) and 1b)
with ramification locus and the components of positive genus 
obtained by blow up.}
\label{1abfibers}
\end{center}
\end{figure}

2)a): We will keep using the polynomials $f(x_0)$ and $g(x_0)$
of the proof of part 1)a).
Let $\AP$ be the smooth $R$-model for $\PK$ corresponding to the 
coordinate $x_0$ and 
$\mathcal X$ its normalization in $K(X)$. 
Then the singularities of $\mathcal X$ are closed points and their image
on the special fiber $\Pk$ of $\AP$ is 
contained in the locus $Z=\{x_0|\bar{g}(x_0)=0\}$. Now
$$\bar{g}(x_0)=x_0^2-x_0\left( \frac{\C+1}{\B+\C+1}\right)^- \quad 
\mbox{\rm and} \quad
Z=\{\bar{d}=(\frac{\C+1}{\B+\C+1})^-,0\}\quad \mbox{\rm with } \bar{d} 
\not=0.$$
We will verify that a suitable blow up on $\AP$ in the point 
$\bar{d}$ will induce a component of genus $(p-1)/2$ in the stable model 
of $X$:
Set $d=\frac{\C+1}{\B+\C+1}$; then this notation is consistent with the above
use of $\bar{d}$. We will blow up the ideal $(x_0-d,b)$ on $\AP$ where
$b=\tau^{1/2}$. With $N$ as before the equation for the generic fiber 
with respect to the coordinate $x_1=b/(x_0-d)$ is \eqref{blow}.

Now observe that $\bar{d} \notin \{0,1\}$, i.e.\ it is not a zero of
$\bar{f}(x_0)$, so $f(d)$ is a unit in $R$. Also by construction $f'(d)=0$
and $f^{(2)}(d)$ is a unit.
Using Propositions \ref{ap} and \ref{stich} we conclude that $x_1=b/(x_0-d)$ 
corresponds to a component
of genus $(p-1)/2$ in the stable reduction of X.
At this point, by standard properties of semi-stable models,
it is already clear that the special fiber of the stable 
reduction is as stated in the theorem. To actually also obtain an equation of 
the second component one passes to a smooth $R$-model for $\PK$ on which
the two points $x_0=0$ and $x_0=\llll$ have distinct specializations by
change of coordinates $x_0 \mapsto x_0/\llll$.
Then one proceeds as before to get a different component of the same genus.

2)b1): Consider the smooth $R$-model for $\PK$ corresponding to
$x_1=x_0/\tau$. We claim that its normalization in $K(X)$ is a smooth model
for $X$. With respect to $x_1$ the branch locus is 
$\{0,\llll/\tau, \infty, 1/\tau\}$. Now the claim follows from 
\cite{L1} Example 1 and Proposition 4. 
Alternatively one can write out the equations
as done before. The statement on the $p$-rank is immediate from the
Deuring-Shafarevich Formula (cf.\ \cite{Cr} Corollary 1.8.)
once one observes that the special fiber of the smooth model is a
Galois cover of $\Pk$ ramified in two points.

\begin{figure}[h]
\begin{center}
\input{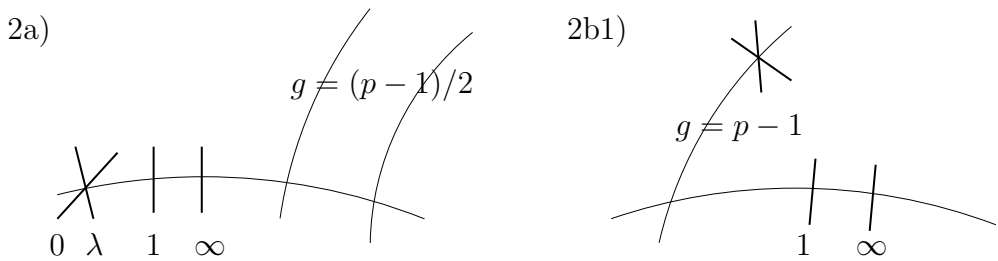}
\caption{Typical picture of the fibers in Theorem \ref{deuring}' 2a) and 2b1)
with ramification locus and the components of positive genus 
obtained by blow up.}
\label{2ab1fibers}
\end{center}
\end{figure}

2)b2): Consider the model for $\PK$ corresponding to $x_1=x_0/\llll^{1/2}$.
With respect to $x_1$ the branch locus is 
$\{0,\llll^{1/2},\infty, 1/\llll^{1/2}\}$.
The fact that the stable model is a Mumford curve now follows exactly
as in \cite{L1} Example 5 (with $\tau^2$ replaced by $\llll^{1/2}$).

2)b3): What follows is the most delicate part of the proof.
To symmetrize the position of the branch points we consider the smooth
$R$-model for $\PK$ corresponding to $x_1=x_0/(x_0+\llll^{1/2})$ and write
$\E=\llll^{1/2}/(1+\llll^{1/2})$.
With respect to the coordinate $x_1$ the branch locus becomes
$\{0,\E,1,1-\E\}$. The equation for the 
generic fiber is
\begin{equation} \label{2}
z_1^p=F(x_1)=x_1(x_1-\E)^{p-1}(x_1-1)^{p-\B}(x_1-1+\E)^{\B}.
\end{equation}
Further we set
\begin{equation}\label{ha} h(x_1)=x_1(x_1-1) \end{equation}

Let ${\mathcal Y}_k$ be the special fiber of the normalization of the new
model (corresponding to $x_1$) in $K(X)$.
Using Proposition \ref{ap} 3) and Remark \ref{app}
one shows that ${\mathcal Y}_k$
is given birationally by the equation
$$y^p=t(x_1)=-\left( \frac{h(x_1)^p-F(x_1)}{\llll^{1/2}} \right)^-.$$ 
Further a direct computation yields
$$t(x_1)=(x_1-1)^{p-1}x_1^{p-1}((\bar{\B}+1)x_1-1) \notin k[x]^p$$
$$t'(x_1)= (x_1-1)^{p-2}x_1^{p-2}(-(\bar{\B}+1)x_1^2+2x_1-1)$$

We conclude that the singularities on $\mathcal Y$ 
map to the locus $\{x_1|-(\bar{\B}+1)x_1^2+2x_1-1=0\}$ if $\bar{\B}+1 \not=0$ 
and, after minor modifications, 
to $\{x_1|2x_1-1=0\} \cup \{\infty\}$ if $\bar{\B}+1=0$.
Notice that these loci are disjoint from $\{0,1\}$, the locus where the
branch points of the generic fiber specialize to.

\begin{figure}[h]
\begin{center}
\input{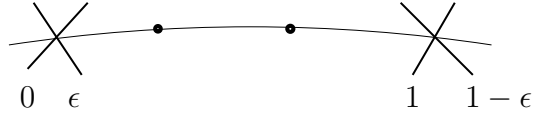}
\caption{The symmetrized model for $\PK$, corresponding to $x_1$, as introduced in the proof 
of case 2b3). Notice that the image of the singular locus, marked by dots, is 
disjoint from the horizontal branching.}
\label{2b3fibers}
\end{center}
\end{figure}

Strictly speaking introducing $x_1=x_0/(x_0+\llll^{1/2})$ is not necessary but it
makes the following constructions less ad hoc. Also we will only consider 
the case that $\bar{\B}+1 \not=0$ the other case can be treated in a similar
manner. At this point we again have to distinguish two cases:

First assume $0 < v(\llll^{1/2}) \le v(p^{(p-2)/(p-1)})$. 
Consider the derivative
$$F'(x_1)=(x_1-\E)^{p-2}(x_1-1)^{p-\B-1}(x_1-1+\E)^{\B-1}$$
$$\left( 2px_1^3-(3p-\E(p-1-\B))x_1^2+x_1(\E^2(\B-1-p)+p+2\E)+\E(\E-1) 
\right)$$


We set $g(x_1)=$
$$\frac{1}{\llll^{1/2}}
\left( 2px_1^3-(3p-\E(p-1-\B))x_1^2+x_1(\E^2(\B-1-p)+p+2\E)+\E(\E-1) \right).$$


Notice that $g(x_1) \in R[x_1]$ and the zeros of $\bar{g}(x_1)$ are contained 
in the image of the singularities of ${\mathcal Y}_k$ above. 
Let $d$ be a zero of $g(x_1)$ not specializing to $\infty$.
There are exactly two of these and they specialize to distinct simple zeros
of $\bar{g}(x_1)$. 
Choose $b$ such that $v(b^2\llll^{1/2})=v(\tau)$ and consider the 
coordinate $x_2=b/(x_1-d)$. Then the equation for the
special fiber becomes
\begin{equation}\label{3}
z_2^p=\sum_{i=0}^{2p}b^i \frac{F^{(i)}(d)}{i!F(d)}x_2^{2p-i}
\end{equation}
Observe that $F(d)$ is a unit, $F'(d)=0$ and $v(F^{(2)}(d))=\llll^{1/2}$ 
because
$d$ was chosen to be a simple zero. Also $v(b^iF^{(i)}(d)/i!)>v(\tau)$ for 
$3 \le i$ and $(p,i)=1$. Further the condition 
$0 < v(\llll^{1/2}) \le v(p^{(p-2)/(p-1)})$ implies that $v(b^p) \ge v(\tau)$.
With $h_2(x_2):=x_2^2$ Proposition \ref{ap}, used as indicated in Remark 
\ref{app},  
shows that this blow up yields an irreducible component
in the stable reduction of $X$. Further, by
Proposition \ref{stich}, the  
genus of this component is $(p-1)/2$.  
Doing the same for the other choice of $d$ will yield the 
second component in the stable model. Having introduced the 
symmetrized model
at the beginning of the proof it is clear that this second component is
different from the one just computed.

Now assume $v(\tau) > v(\llll^{1/2}) \ge v(p^{(p-2)/(p-1)})$. 
With $x_1$ as in equation \eqref{2} and $h(x_1)$ as defined in \eqref{ha}
write $$T_0(x_1)=\frac{h(x_1)^p-F(x_1)}{\llll^{1/2}}$$
We have seen above that $\bar{T}_0'(x_1)=-t'(x_1)$ has two distinct zeros 
outside
$\{0,1\}$. Now let $d$ be one of the two zeros of $T_0'(x_1)$ not specializing 
to $\{0,1\}$. Choose $b$ such that $v(b^2\llll^{1/2})=v(\tau)$.
In what follows we will show that 
the blow up corresponding to $x_2=b/(x_1-d)$ will yield a component of the 
stable model. As before the generic fiber is given birationally by 
$$z_2^p=\tilde{F}(\U)={\U}^{2p}+b\frac{F'(d)}{F(d)}\U^{2p-1}
+b^2\frac{F^{(2)}(d)}{2F(d)}\U^{2p-2} 
+ \dots + b^p\frac{F^{(p)}(d)}{p!F(d)}\U^p+ \dots$$
Set $$\tilde{h}(\U)=\U^2+b\frac{h'(d)}{h(d)}\U+\frac{b^2}{h(d)}$$
Notice that 
$$\tilde{h}(\U)^p=\U^{2p}+b\frac{(h(x_1)^p)'|_{x_1=d}}{h(d)^p}\U^{2p-1}
+b^2\frac{(h(x_1)^p)^{(2)}|_{x_1=d}}{2!h(d)^p}\U^{2p-2}+ \dots $$
We claim that $v(\tilde{h}(\U)^p-\tilde{F}(\U))=v(\tau)$.
With $T(x_1)=h(x_1)^p-F(x_1)$ we have $T'(x_1)=ph(x_1)^{p-1}h'(x_1)-F'(x_1)$.
Further $v(T'(x_1))=v(\llll^{1/2})$ because $\bar{T}_0'(x_1)=-t'(x_1)\not\equiv 0$.
Certainly $v(F'(x_1)) \le v(p)$, so
we get that $v(F'(x_1))=\mbox{\rm min}\{v(\llll^{1/2}),v(p)\}$.
 
Also $v(T(x_1))=v(\llll^{1/2})$ therefore $h(d)^p=F(d)+r$ with 
$v(r) \ge v(\llll^{1/2})$.
Now we are ready to analyze $\tilde{h}(\U)^p-\tilde{F}(\U)$.
Consider the coefficient of $\U^{2p-i}$:
$$b^i\frac{(h(x_1)^p)^{(i)}|_{x_1=d}}{i!h(d)^p}-b^i\frac{F^{(i)}(d)}{i!F(d)}$$
Multiplying this with the unit $h(d)^p$ we get
$$\frac{b^i}{i!}\left( (h(x_1)^p)^{(i)}|_{x_1=d}
-F^{(i)}(d)\frac{h(d)^p}{F(d)} \right)$$
$$=\frac{b^i}{i!}\left( (h(x_1)^p-F(x_1))^{(i)}|_{x_1=d}\right) 
+b^ir\frac{F^{(i)}(d)}{i!F(d)}.$$
Now by construction $(h(x_1)^p-F(x_1))'|_{x_1=d}=0$ and $v(brF'(d)) > v(\tau)$.
(It is this last inequality that only is valid if we put a lower bound 
on $v(\llll^{1/2})$).
Further the coefficient of $\U^{2p-2}$ has order exactly $v(\tau)$
because $\bar{d}$ is a simple zero of $\bar{T}_0'(x_1)$. 
Finally the coefficients of
all other $\U$-powers have orders $>v(\tau)$. 
Again Proposition \ref{ap} case 2 yields
a component of genus $(p-1)/2$. The other component is obtained the same way
using the second zero of $T_0'(x_1)$ mentioned above. \ed \\

Next we treat the case of $p=3$. Otherwise
we keep the assumptions of the previous theorem. Here the result will look
slightly different.

\begin{theorem} \label{p3}
There are three possibilities for the special fiber $X_k$ of the stable model:

\begin{itemize}
\item[1)] $X_k$ consists of a single irreducible component, i.e.\ 
$X$ has potentially good reduction iff $\C+1=p$ and $v(\llll)=v(\tau^2)$.
In this case $X_k$ has $p$-rank $p-1$.

\item[2)] $X_k$ has exactly two irreducible components, both of genus zero and
intersecting each other in $p$ distinct points iff $\C+1=p$ and 
$v(\llll)>v(\tau^2)$. In this case $X$ is a Mumford curve.

\item[3)] In all remaining cases $X_k$ has two irreducible components, each
of genus $(p-1)/2$, intersecting in precisely one point.
\end{itemize}
\end{theorem}

\noindent{\it Proof.} We follow the same path as for $p>3$ distinguishing
the cases $v(\llll)=0$ and $v(\llll) \ge 0$.
The first case is a consequence of \cite{M} and the algorithm given there
also will yield equations. Observe that it is a priory clear that potentially 
good reduction with $p$-rank zero is not an option, using genus formulas.

In the case of $v(\llll) \ge 0$ the previous proof carries over word by
word. \ed\\

We would like to compare the above to a classic result of Deuring 
on elliptic curves which we first recall:\\

\noindent{\bf Theorem (Deuring)}{\it
Let $X/K$ be an elliptic curve given by $y^2=x(x-1)(x-\llll)$
with $j$-invariant
\[j(E)=2^8\frac{(\llll^2-\llll+1)^3}{\llll^2(\llll-1)^2} .\] 
Then $X$ has potentially good reduction if and only if $j(E)\in R$.}\\

To see why Theorem \ref{deuring} is an analog of this consider 
Deuring's criterion in
the case that the residue characteristic is $p=2$. This is the only case
which involves wild ramification. Then it is a statement
about reduction of cyclic covers of degree $p=2$, ramified at four points.
Deuring's result is not a special case of Theorem \ref{deuring} but can be
proved using methods similar to those used above.
Note that at the beginning we defined 
$j(\llll)$ for $z_0^p=x_0(x_0-1)^{\B}(x_0-\llll)^{\C}$ and in the special case of
$\B=\C=1$ this yields $j(\llll)=4p^{-2p/3(p-1)}(\llll^2-\llll +1)$ which
is quite similar to $j(E)$ for elliptic curves.

\begin{remark} \rm 
In Theorem \ref{deuring} it would be desirable to have a set of invariants
such as $j(\llll)$ which determine the stable reduction of $X$ through
their absolute values. This has been achieved
by Q.Liu for genus $2$ curves (cf.\ \cite{Li} Theorem 1). 
In our setup it seems 
not obvious how to find such a set of invariants if possible at all.
A reason for this is that here there is no modular interpretation in the
sense that not every curve of genus $(m-1)(p-1)/2$ arises as a $p$-cyclic cover
of the affine line.
\end{remark}

\section{An Application} \label{exo}

The following application of Theorem \ref{deuring} was pointed out by David
Harbater. Besides being of interest on its own, it also has significance
in the theory of lifting Galois covers from positive characteristic to
characteristic zero.

\begin{proposition} \label{qwerty} For $p>3$
consider the cover $X \la \Pk$ given birationally by 
\begin{equation} \label{exx}
z^p=(x-c_1)^{p-1}(x+c_1)(x-c_2)^{p-1}(x+c_2)
\end{equation}
where $c_1,-c_1,c_2,-c_2 \in R$ are all distinct.
Then this cover can not have potentially good reduction of type 1a.
\end{proposition}

\begin{proof} We first bring equation \eqref{exx} in the shape
required to apply Theorem \ref{deuring}. Therefore we introduce the coordinate
$$x_0=\frac{c_1-c_2}{2c_1} \frac{x+c_1}{x-c_2}.$$
With respect to $x_0$ the generic fiber is given by
\begin{equation}
z_0^p=x_0(x_0-1)^{p-1}(x_0+\frac{(c_1-c_2)^2}{4c_1c_2})
\end{equation}

In what follows we will assume that the cover has potentially good reduction
of type 1a) and produce a contradiction. 
Observe that Theorem \ref{deuring}' 1)a) implies that there exists a smooth $R$-model
for $\PK$ such that the points of the branch locus of $X \la \PK$ specialize
to distinct points on the closed fiber.
Now, with respect to $x_0$, the branch locus contains $\{0,1,\infty\}$ and hence
the fourth point can not specialize to the set $\{0,1,\infty\}$ on the closed fiber.
We write $\llll=-\frac{(c_1-c_2)^2}{4c_1c_2}$ and conclude $\llll \in R-\{0,1\}$ and
$\bar{\llll}\not=1$. Now we can apply Theorem \ref{deuring} 1)a) to obtain 
$v(j(\llll)) \ge 0$. So
$j(\llll)=p^{-2p/(3(p-1))}(\llll^2(\B +1)^2-2\llll
(\B+\C+1-\B\C)+(\C+1)^2)$ and with
$\B=p-1$ and $\C=1$ we get
$j(\llll)= p^{-2p/(3(p-1))}(\llll^2p^2-4\llll+4)$
This yields $v(\llll-1) \ge 0$ - a contradiction, so this case
doesn't occur.

\begin{remark} \rm In the situation of Proposition \ref{qwerty} potentially
good reduction of type 1b is possible. One easily checks that it is the case 
for $c_1=1$ and $c_2=-\tau-1$ or $c_1=1$ and $c_2=\tau^2$.
\end{remark}

\end{proof}

\bibliographystyle{amsplain}

\begin{thebibliography}{AAAAA}





\bibitem[Cr]{Cr} Crew, R.M.: Etale $p$-covers in characteristic $p$.
Compositio Math.{\bf 52}, 31-45 (1984).

\bibitem[DM]{DM}
Deligne, P., Mumford, D.: The irreducibility of the
space of curves of given genus. Inst. Hautes Etudes Sci. Publ. Math. 
{\bf 36}, 75-109 (1969).












\bibitem[L1]{L1}
Lehr, C.: Reduction of $p$-cyclic Covers of the Projective Line.
Manuscripta Math. {\bf 106}, 151-175 (2001).

\bibitem[L2]{L2} Lehr, C. {Effective methods for vanishing cycles of 
$p$-cyclic covers of the $p$-adic line}. 
J. Algebra {\bf 271}, no. 1, 407-425 (2004).

\bibitem[Li]{Li}
Liu, Q.: Courbes stables de genre $2$ et leur sch\'ema de modules.
Math.Ann. {\bf 295}, 201-222 (1993).



\bibitem[M]{M}
Matignon, M.: Vers un algorithme pour la r\'eduction stable des rev\^etements
$p$-cycliques de la droite projective sur un corps $p$-adique.
Math. Ann. {\bf 325}, no. 2, 323-354 (2003).

\bibitem[R]{R1}
Raynaud, M.: $p$-groupes et r\'eduction semi-stable des courbes. The
Grothendieck Festschrift, Vol.3, Basel-Boston-Berlin: Birkh\"auser 1990.



 

 
\bibitem[St]{St} Stichtenoth, H.: Algebraic Function Fields and
Codes. Berlin-Heidelberg-New York: Springer 1993. 



\bibitem[Za]{Za} Zapponi, L.: 
Specialization of polynomial covers of prime degree. 
Pacific J. Math. {\bf 214}, no. 1, 161-183 (2004).



\end{thebibliography}

\end{document}